\documentclass[11pt]{article} 

\usepackage{latexsym} 
\usepackage{amsmath}
\usepackage{amssymb} 
\usepackage{amscd} 
\usepackage{graphicx}
\begin{document} 
\newtheorem{Th}{Theorem}[section]
\newtheorem{Cor}{Corollary}[section]
\newtheorem{Prop}{Proposition}[section]
\newtheorem{Lem}{Lemma}[section]
\newtheorem{Def}{Definition}[section]
\newtheorem{Rem}{Remark}[section]
\newtheorem{Ex}{Example}[section]
\newtheorem{stw}{Proposition}[section]


\newcommand{\bet}{\begin{Th}}
\newcommand{\ent}{\stepcounter{Cor}
   \stepcounter{Prop}\stepcounter{Lem}\stepcounter{Def}
   \stepcounter{Rem}\stepcounter{Ex}\end{Th}}


\newcommand{\bec}{\begin{Cor}}
\newcommand{\enc}{\stepcounter{Th}
   \stepcounter{Prop}\stepcounter{Lem}\stepcounter{Def}
   \stepcounter{Rem}\stepcounter{Ex}\end{Cor}}
\newcommand{\bep}{\begin{Prop}}
\newcommand{\enp}{\stepcounter{Th}
   \stepcounter{Cor}\stepcounter{Lem}\stepcounter{Def}
   \stepcounter{Rem}\stepcounter{Ex}\end{Prop}}
\newcommand{\bel}{\begin{Lem}}
\newcommand{\enl}{\stepcounter{Th}
   \stepcounter{Cor}\stepcounter{Prop}\stepcounter{Def}
   \stepcounter{Rem}\stepcounter{Ex}\end{Lem}}
\newcommand{\bef}{\begin{Def}}
\newcommand{\enf}{\stepcounter{Th}
   \stepcounter{Cor}\stepcounter{Prop}\stepcounter{Lem}
   \stepcounter{Rem}\stepcounter{Ex}\end{Def}}
\newcommand{\ber}{\begin{Rem}}
\newcommand{\enr}{
   \stepcounter{Th}\stepcounter{Cor}\stepcounter{Prop}
   \stepcounter{Lem}\stepcounter{Def}\stepcounter{Ex}\end{Rem}}
\newcommand{\bee}{\begin{Ex}}
\newcommand{\ene}{
   \stepcounter{Th}\stepcounter{Cor}\stepcounter{Prop}
   \stepcounter{Lem}\stepcounter{Def}\stepcounter{Rem}\end{Ex}}
\newcommand{\Proof}{\noindent{\it Proof\,}:\ }

\title{Singular curves of hyperbolic $(4, 7)$-distributions of type $C_3$}

\author{Goo Ishikawa and Yoshinori Machida} 

\date{ }

\maketitle

\renewcommand{\thefootnote}{\fnsymbol{footnote}}


\newcommand{\R}{\mathbf{R}}
\newcommand{\C}{\mathbf{C}}
\newcommand{\QED}{\hfill $\Box$}


\begin{abstract}
A distribution of rank $4$ on a $7$-dimensional manifold is called a $(4, 7)$-distribution 
if its local sections generate the whole tangent space by taking Lie brackets once. 
Singular curves of $(4, 7)$-distributions are studied in this paper. 
In particular the class of hyperbolic $(4, 7)$-distributions of type $C_3$ is introduced and 
singular curves are completely described via prolongations for them.

\end{abstract}

\maketitle

\section{Introduction}
\label{Introduction}

Let $X$ be a differentiable (= $C^\infty$) manifold of dimension $7$ and $D \subset TX$
a distribution (i.e., subbundle) of rank $4$. 
The distribution $D$ is called a {\it $(4, 7)$-distribution} if the local sections to $D$ generate 
the whole tangent bundle $TX$ by taking Lie bracket once, i.e., if ${\mathcal {TX}} = [{\mathcal D}, {\mathcal D}]$. 
See \S \ref{(4,7)-distributions} for details.   

Richard Montgomery studies $(4, 7)$-distributions in \cite{Montgomery}. 
He defines elliptic and hyperbolic $(4, 7)$-distributions 
and shows that, to any elliptic (resp. hyperbolic) $(4, 7)$-distribution $D \subset TX$, there corresponds uniquely a conformal definite (resp. $(2, 2)$) metric on $D$. Moreover he solves the classification problem of such kinds of distributions 
in the sense of Cartan-Tanaka theory (\cite{Montgomery} \S 6.10 and \S 7.12). 

We study in this paper singular curves of $(4, 7)$-distributions. 
Singular curves provide a characteristic information on distributions and the class of singular curves is 
regarded as an important invariant of the isomorphism classes of distributions. 
For singular curves see \S \ref{Singular curves}. 

The elliptic $(4, 7)$-distribution has no non-trivial singular curves (see \cite{Montgomery} \S 5.6 
and \S \ref{(4,7)-distributions} of the present paper). Thus we are led naturally to 
ask for singular curves of hyperbolic $(4, 7)$-distributions. 

Given a distribution $D \subset TX$, 
the velocity vectors of $D$-singular curves form 
a cone field over $X$, which 
we call the {\it singular velocity cone} of $D$ and denoted by ${\mbox {\it{SVC}}} \, (\subset TX)$. 
Any direction $[u]$ in $T_xX$, where $(x, u) \in {\mbox {\it{SVC}}}$, is called a {\it singular direction} of $D$ 
(Definition \ref{singular-velocity-cone}). 
Then we show there exists a $(2, 2)$-metric on $D$ with the null-cone field $C \subset TX$ such that 
${\mbox {\it{SVC}}} \subseteq C$ and therefore any $D$-singular curve is necessarily an integral curve of $C$ 
(Proposition \ref{necessary-condition}). 
The $(2, 2)$ metric is exactly the metric found in \cite{Montgomery}. 
We call the cone field $C \subset TX$ the {\it characteristic cone} of $D \subset TX$. 
Any direction $[u] \in P(C_x)$ is called 
a {\it characteristic direction} for $(x, u) \in C, u\not= 0$. Note that, outside of the zero section, 
$C$ is a submanifold of $TX$ of dimension $10$. 

Our fundamental problem is to ask whether 
a characteristic direction in $C$ is a singular direction or not,  
and to determine ${\mbox {\it{SVC}}}$ for $(4, 7)$-distributions. 

To approach to these problem, we prolong the cone field $C \subset D$ which 
is intrinsically defined by $D$. 
Then we get a distribution $E$ of rank $3$ 
on the $9$-dimensional manifold $Z = PC = (C \setminus (\mbox{\rm zero section}))/\R^\times$, 
the set consisting of all directions of $C$, with the canonical projection $\pi : Z \to X$, $(x, [u]) \mapsto x$. 
Then, in our case,  we see that the small growth of $E$ at each point of $Z$ turns out to be 
$(3, 5, 7, 8, 9)$ or $(3, 5, 7, 9)$, depending on a point of $Z$ 
(see \S \ref{Prolongations of hyperbolic $(4, 7)$-distributions}). 

\bef
\label{C_3-type-def}
{\rm 
Let $x \in X$ and $u \in C_x \setminus \{ 0\}$. 
Then we call the direction $[u] \in P(C_{x})$ {\it of type $C_3$} if 
the small growth of $E$  is equal to 
$(3, 5, 7, 8, 9)$ on a neighbourhood of $(x, [u])$ in $Z =PC$. 

We call a $(4, 7)$-distribution $D \subset TX$ {\it of type $C_3$} 
if for any $(x, [u]) \in Z$, $[u]$ is of type $C_3$. 
}
\enf

A typical example of $(4, 7)$ distribution of type $C_3$ is given by 
the canonical distribution $D \subset TX$ on a Grassmannian $X$ which consists of isotropic 
$2$-dimensional subspaces in the $6$-dimensional symplectic vector space $(\R^6, \Omega)$. 
Note that $X$ is a homogeneous space of Lie group of type $C_3$ (See \S\ref{Prolongations of hyperbolic $(4, 7)$-distributions} and \S\ref{Isotropic-Grassmannian}). 

We call a $C^\infty$ immersive $D$-singular curve a {\it $D$-singular path}. 
In this paper we show 

\bet
\label{main-theorem} 
Let $D \subset TX$ be a hyperbolic $(4, 7)$-distribution. Then 

{\rm (1)} If the direction $[u] \in P(C_{x})$ is of type $C_3$, then 
there exists a $D$-singular path $\gamma : (\R, 0) \to X$ with $\gamma(0) = x, [\gamma'(0)] = [u]$. 
In particular $u \in {\mbox {\it{SVC}}}$. 

{\rm (2)} If $D$ is of type $C_3$, then for any $[u] \in P(C_{x})$, there exists 
a $D$-singular path $\gamma : (\R, 0) \to X$ with $\gamma(0) = x, [\gamma'(0)] = [u]$. 
In particular ${\mbox {\it{SVC}}} = C$. 
\ent

Montgomery shows also that an orientation of $D$ is uniquely determined from $D$. The orientation is determined 
depending on that the three plane of signature $(2, 1)$ in $D^*\wedge D^*$ defined by the dual curvature map is self-dual or anti-self-dual (\cite{Montgomery}), which is read also from singular curves of $D$ in our case. 

To show Theorem \ref{main-theorem}, we analyse $E$-singular curves of the prolongation $(Z, E)$ of an intrinsic $3$-dimensional cone $C \subset D$ defined from $(X, D)$. 
Then we observe, as a characteristic property of a hyperbolic $(4, 7)$-distribution $(X, D)$ of type $C_3$, 
that there exists canonically a sub-distribution $F \subset E$ of rank $2$ on $Z$ which characterises $D$-singular curves. 

We characterise $D$-singular curves completely for a hyperbolic $(4, 7)$-distribution $(X, D)$ of 
type $C_3$. 

\bet
\label{main-theorem2} 
Let $(X, D)$ be a hyperbolic $(4, 7)$-distribution of type $C_3$ and $(Z, E)$ the prolongation by 
the characteristic cone field $C \subset D$ with the projection $\pi : Z \to X$. 
Then there exists a subbundle $F \subset E$ of rank $2$ such that 

{\rm (1)} Any $F$-integral curve $c : I \to Z$ with $(\pi\circ c)'(t) \not= 0$, 
for its projection $\pi\circ c : I \to X$ by $\pi : Z \to X$, is an $E$-singular path and 
$\pi\circ c$ is a $D$-singular path. 

{\rm (2)} Let $(x, [u]) \in Z$. Then any $D$-singular path $\gamma : (\R, 0) \to X$ 
with $\gamma(0) = x$ and $[\gamma'(0)] = [u]$ 
is obtained as the projection $\pi\circ\widetilde{\gamma}$ 
of an $F$-integral curve $\widetilde{\gamma} : (\R, 0) \to L$
with $\widetilde{\gamma}(0) = (x, [u])$. 
\ent

By Theorem \ref{main-theorem2}, we see, in Theorem \ref{main-theorem}, 
that a $D$-singular curve with given an initial point and an 
initial direction is not uniquely determined. 

In the typical case of hyperbolic $(4, 7)$-distributions of type $C_3$ (\S\ref{Isotropic-Grassmannian}), 
the distribution $F$ appeared in Theorem \ref{main-theorem2} induces an {\it Engel foliation} on $Z$, 
i.e., a foliation ${\mathcal F}$ on $Z$ with four dimensional leaves such that $F\vert_L \subset TL$ 
and $F\vert_L$ has small growth $(2, 3, 4)$ for any leaf $L$ of ${\mathcal F}$ 
(see Remark \ref{Engel-foliation} of \S\ref{Isotropic-Grassmannian}). 
Recall that a distribution of rank $2$ on a $4$-dimensional manifold $L$ is called an {\it Engel distribution} 
if its small growth is $(2, 3, 4)$ (\cite{Montgomery}). 

\

In the following sections \S  \ref{Singular curves}, \S \ref{(4,7)-distributions} 
and \S \ref{Prolongations of hyperbolic $(4, 7)$-distributions}, 
we review on singular curves and on $(4, 7)$-distribution and related topics. 
In particular we study the prolongation $(Z, E)$ along $C \subset TX$ 
for a hyperbolic $(4, 7)$-distributions $(X, D)$.
We prove Theorem \ref{main-theorem} and Theorem \ref{main-theorem2} in \S \ref{Proof of Theorem}. 
In the last section \S \ref{Isotropic-Grassmannian}, we give the model of hyperbolic $(4, 7)$-distribution $D$ on the isotropic Grassmannian on $6$-dimensional symplectic space (\cite{Montgomery0}). Moreover we study its
perturbations as an example of hyperbolic $(4, 7)$-distributions of mixed type. 

\

For geometric studies on singular curves related to several simple Lie algebras, 
see \cite{Tanaka, SY} and \S \ref{(4,7)-distributions} of this paper. 

Montgomery \cite{Montgomery0} conjectured 
and Kry{\' n}ski \cite{Krynski} showed that generic distributions of corank at most $3$ 
are determined by their singular curves except for several exceptional cases. 
The case of $(4, 7)$-distributions is one of exceptions, because of the existence of elliptic $(4, 7)$-distributions. 
Then it is an interesting problem to study singular curves of generic hyperbolic $(4, 7)$-distributions, which is not treated in this paper. 

For complex analytic distributions, we can define the notion of singular curves in terms of constrained Hamiltonian equations 
as explained in \S \ref{Singular curves}, not in terms of control theory. Then it can be proved a similar result to 
Theorems \ref{main-theorem}, \ref{main-theorem2} 
on singular velocity cones and singular curves of holomorphic \lq\lq non-degenerate\rq\rq $(4, 7)$-distributions 
of type $C_3$ which have prolongations with growth $(3, 5, 7, 8, 9)$. 

All manifolds, mappings, vector fields, distributions, etc. are assumed of class $C^\infty$ 
unless otherwise stated. 

\section{Singular curves of distributions}
\label{Singular curves}

Let $D \subset TX$ be any distribution on a manifold $X$ and $x_0 \in X$. 
Let us consider the space ${\mathcal C}_{x_0}$ of $D$-integral curves, 
which consist of all absolutely continuous curves 
$\gamma : [a, b] \to X$ with $\gamma(a) = x_0, \gamma'(t) \in D_{\gamma(t)}$ for almost every $t \in [a, b]$. 
A curve $\gamma : [a, b] \to X$ starting at $x_0$ is called 
a {\it singular curve} of $D$ or a {\it $D$-singular curve}, 
if it is a critical point of the endpoint mapping 
${\mbox{\rm end}} : {\mathcal C}_{x_0} \to X$ defined by ${\mbox{\rm end}}(\gamma) = \gamma(b)$. 
Therefore a singular curve is a curve such that we can't control its endpoint infinitesimally. 
Singular curves are characterised locally by a constrained Hamiltonian equation. 

Let $\xi_1, \dots, \xi_r$ be a local frame of $D$. Just for simplicity, suppose $\xi_1, \dots, \xi_r$ are defied on $X$. 
Then we define the Hamiltonian function on the cotangent bundle $T^*X$ of $D$ by 
$$
H(x, p; u) := u_1 H_{\xi_1}(x, p) + \cdots + u_r H_{\xi_r}(x, p), 
$$
where $H_{\xi_i}(x, p) = \langle p, \xi_i(x)\rangle$ for $(x, p) \in T^*X$ and $u_1, \dots, u_r$ are control parameters. 
Then consider the differential equation on $(x, p, u) = (x(t), p(t), u(t))$, 
$$
\dot{x} = \frac{\partial H}{\partial p}(x, p), \quad \ \dot{p} = - \frac{\partial H}{\partial x}(x, p)  \quad \cdots\cdots\cdots (*)
$$
with constraints $H_{\xi_1}(x, p) = 0, \dots, H_{\xi_r}(x, p) = 0$. 
The constraints means that $p(t) \in D_{x(t)}^\perp$. 
We say that $(x(t), p(t))$ is a {\it singular bi-characteristic} if there exist $u_1(t), \dots, u_r(t)$ depending on $t$ 
such that $(x(t), p(t), u(t))$ satisfies the equation (*) with the constraints. 
A curve $\gamma(t) = x(t)$ is a $D$-singular curve if and only if there exist a singular bi-characteristic 
$\Gamma(t) = (x(t), p(t))$ covering $x(t)$ with $p(t) \not= 0$ (see \cite{AS, LS}). 
Note that a $D$-singular curve satisfies 
\[
\dot{x} = u_1(t)\xi_1(x(t)) + \cdots + u_r(t)\xi_r(x(t)). 
\]
with the controls $u_1(t), \dots, u_r(u)$. 

Later we use the following notation and fact: 
For a vector field $\xi$ over $X$ we denote by $\overrightarrow{H_{\xi}}$ the Hamiltonian vector field over $T^*X$ with 
Hamiltonian $H_{\xi}$. Note that the vector field $\overrightarrow{H_{\xi}}$ over $T^*X$ covers the vector field $\xi$ 
over $X$ via the canonical projection $T^*X \to X$. 
If $(x(t), p(t), u(t))$ satisfies the equation (*), 
then we have along the curve 
\[
\frac{d}{dt}H_{\xi}(x(t), p(t)) = u_1(t) H_{[\xi_1, \xi]}(x(t), p(t)) + \cdots + u_r(t) H_{[\xi_r, \xi]}(x(t), p(t)), 
\]
for any vector field $\xi$ over $X$ (see \cite{AS, LS}). 

Now we define the key notion in this paper: 

\bef
\label{singular-velocity-cone}
{\rm 
The {\it singular velocity cone} ${\mbox{\it{SVC}}} \subset TX$ is defined as the set of 
$(x, u) \in TX$ such that there exists a $C^\infty$ singular curve $\gamma : (\R, 0) \to X$ such that 
$\gamma(0) = x, \gamma'(0) = u$. 
}
\enf

Note that ${\mbox{\it{SVC}}} \subset D$ and ${\mbox{\it{SVC}}}$ is a cone, i.e., it is invariant under the $\R^\times$-action on $TX$. In fact, if $(x, u) \in TX$ is realised by a singular curve $\gamma(t)$ 
at $t = 0$ and $c \in \R^\times$, then $(x, cu)$ is realised by the singular curve $\gamma(ct)$.

\section{$(4, 7)$-distributions and their singular curves}
\label{(4,7)-distributions}

Let $X$ be a $C^\infty$ manifold of dimension $7$ and $D \subset TX$ a $C^\infty$ distribution (= subbundle of $TX$) 
of rank $4$. 

We call $(X, D)$ a {\it $(4, 7)$ distribution} if Lie brackets of local sections to $D$ generate $TX$, 
i.e., if 
$D^{(2)} := [D, D] (= D + [D, D]) = TX$, 
or equivalently, 
if, for any $x_0 \in X$, there exists a local frame $\xi_1, \xi_2, \xi_3, \xi_4$ 
of $D$ in a neighbourhood $U$ of $x_0$ such that 
$\xi_i(x_0)$ and $[\xi_j, \xi_k](x_0)$ generate $T_{x_0}X$. 

According to \cite{Montgomery} we define two classes of 
$(4, 7)$-distributions. 

\

A $(4, 7)$-distribution $(X, D)$ is called {\it elliptic} if there exists a local frame 
$\xi_1, \xi_2, \xi_3, \xi_4$ of $D$ such that 
\[
[\xi_1, \xi_2] \equiv - [\xi_3, \xi_4], \ [\xi_1, \xi_3] \equiv [\xi_2, \xi_4], \ 
[\xi_1, \xi_4] \equiv - [\xi_2, \xi_3], \ 
{\mbox {\rm mod}}. D, 
\]
and that 
\[
\xi_1, \xi_2, \xi_3, \xi_4, \ \xi_5 := [\xi_1, \xi_2], \ \xi_6 := [\xi_1, \xi_3], 
\ \xi_7 := [\xi_1, \xi_4]
\]
form a local frame of $TX$. 

See \cite{Montgomery} for an example of elliptic $(4, 7)$-distribution. 

\

A $(4, 7)$-distribution $(X, D)$ is called {\it hyperbolic} 
if there exists a local frame 
$\xi_1, \xi_2, \xi_3, \xi_4$ of $D$ such that 
\[
[\xi_1, \xi_2] \equiv 0,  \ [\xi_3, \xi_4] \equiv 0, \ [\xi_1, \xi_4] \equiv [\xi_2, \xi_3], 
{\mbox {\rm mod}}. D
\] 
and that 
\[
\xi_1, \xi_2, \xi_3, \xi_4, \ \xi_5 := [\xi_1, \xi_3], \ \xi_6 := [\xi_1, \xi_4], 
\ \xi_7 := [\xi_2, \xi_4]
\] 
form a local frame of $TX$. 

We call such a frame $\xi_1, \xi_2, \xi_3, \xi_4$ an {\it adapted frame}. 
See \S \ref{Isotropic-Grassmannian}, for examples of hyperbolic $(4, 7)$-distributions. 

\

Recall the Hamiltonian equation for $D$-singular curves, 
\[
\left\{ 
\begin{array}{l}
{\displaystyle 
\dot{x} = \ u_1\frac{\partial H_{\xi_1}}{\partial p} + u_2\frac{\partial H_{\xi_2}}{\partial p} + u_3\frac{\partial H_{\xi_3}}{\partial p} + 
u_4\frac{\partial H_{\xi_4}}{\partial p}, 
}
\vspace{0.1truecm}
\\
{\displaystyle 
\dot{p} = - \left(u_1\frac{\partial H_{\xi_1}}{\partial x} + u_2\frac{\partial H_{\xi_2}}{\partial x} + u_3\frac{\partial H_{\xi_3}}{\partial x} + 
u_4\frac{\partial H_{\xi_4}}{\partial x}\right), 
}
\end{array}
\right.
\]
with constraints $H_{\xi_1} = 0, H_{\xi_2} = 0, H_{\xi_3} = 0, H_{\xi_4} = 0$ and $p \not= 0$. 
Note that $\dot{x} = u_1\xi_1(x) + u_2\xi_2(x) + u_3\xi_3(x) + u_4\xi_4(x)$. 
Along the singular bi-characteristic, we have, for $i = 1, 2, 3, 4$, 
$$
0 = \frac{d}{dt}H_{\xi_i} = - \left( u_1 H_{[\xi_i, \xi_1]} + u_2 H_{[\xi_i, \xi_2]} + u_3 H_{[\xi_i, \xi_3]} + u_4 H_{[\xi_i, \xi_4]}\right). 
$$
Thus, if we set
\[
A(t) = \left(
\begin{array}{cccc}
0 & H_{[\xi_1, \xi_2]}  & H_{[\xi_1, \xi_3]} & H_{[\xi_1, \xi_4]}
\\
H_{[\xi_2, \xi_1]}  & 0 &  H_{[\xi_2, \xi_3]} & H_{[\xi_2, \xi_4]}
\\
H_{[\xi_3, \xi_1]}  & H_{[\xi_3, \xi_2]} & 0 &  H_{[\xi_3, \xi_4]}
\\
H_{[\xi_4, \xi_1]}  & H_{[\xi_4, \xi_2]} &  H_{[\xi_4, \xi_3]} & 0 
\end{array}
\right), 
\]
we have necessarily $A(t)u = 0$, if there is a singular curve $x(t)$ satisfying the equation 
with $p(t) \in D^\perp_{x(t)} \setminus \{ 0\}$. 
Here $D^\perp : = \{ (x, p) \in T^*X \mid \langle p, D\rangle = 0 \} \subset T^*X$, the annihilator of $D$. 

In the elliptic cases, the necessary condition becomes 
\[
\left(
\begin{array}{cccc}
0 & H_{\xi_5}  & H_{\xi_6} & H_{\xi_7}
\\
- H_{\xi_5}  & 0 &  -H_{\xi_7} & H_{\xi_6}
\\
- H_{\xi_6}  & H_{\xi_7} & 0 &  - H_{\xi_5}
\\
- H_{\xi_7}  & - H_{\xi_6} &  H_{\xi_7} & 0 
\end{array}
\right)
\left(
\begin{array}{c}
u_1
\\
u_2
\\
u_3
\\
u_4
\end{array}
\right)
= 
\left(
\begin{array}{c}
0
\\
0
\\
0
\\
0
\end{array}
\right), 
\]
equivalently, 
\[
\left(
\begin{array}{ccc}
u_2 & u_3 & u_4 
\\
-u_1 & u_4 & - u_3 
\\
-u_4 & - u_1 & u_2
\\
u_3 & -u_2 & -u_1
\end{array}
\right)
\left(
\begin{array}{c}
H_{\xi_5}
\\
H_{\xi_6}
\\
H_{\xi_7}
\end{array}
\right)
= 
\left(
\begin{array}{c}
0
\\
0
\\
0
\\
0
\end{array}
\right), 
\]
for some $H_{\xi_5}, H_{\xi_6}, H_{\xi_7}$.  
Since $p\not= 0$ and $\xi_1, \dots, \xi_7$ generate $TX$, 
we see $H_{\xi_5}, H_{\xi_6}, H_{\xi_7}$ are not all zero. 
Therefore necessarily we have that $u_i(u_1^2 + u_2^2 + u_3^2 + u_4^2) = 0, i = 1,2,3,4$. 
Since we are working on the real numbers, we have $u_1 = u_2 = u_3 = u_4 = 0$.

\bep {\rm (\cite{Montgomery}) } 
For any elliptic $(4, 7)$-distribution $(X, D)$, 
the singular velocity cone consists of zero vectors : ${\mbox {\it{SVC}}} = {\mbox{\rm (zero section)}}$. 
Therefore there exists no non-constant $D$-singular curve on $X$ in elliptic cases. 
\enp

The necessary condition for the existence of $D$-singular curves on a hyperbolic $(4, 7)$-distribution $(X, D)$ 
reads 
\[
\left(
\begin{array}{cccc}
0 & 0 & H_{\xi_5} & H_{\xi_6}
\\
0  & 0 & H_{\xi_6} & H_{\xi_7}
\\
- H_{\xi_5}  & -H_{\xi_6} & 0 &  0
\\
- H_{\xi_6}  & - H_{\xi_7} & 0 & 0 
\end{array}
\right)
\left(
\begin{array}{c}
u_1
\\
u_2
\\
u_3
\\
u_4
\end{array}
\right)
= 
\left(
\begin{array}{c}
0
\\
0
\\
0
\\
0
\end{array}
\right), 
\]
equivalently, 
\[
\left(
\begin{array}{ccc}
u_1 & u_2 & 0
\\
u_3 & u_4 & 0
\\
0 & u_1 & u_2
\\
0 & u_3 & u_4
\end{array}
\right)
\left(
\begin{array}{c}
H_{\xi_5}
\\
H_{\xi_6}
\\
H_{\xi_7}
\end{array}
\right)
= 
\left(
\begin{array}{c}
0
\\
0
\\
0
\\
0
\end{array}
\right)
\]
for some $H_{\xi_5}, H_{\xi_6}, H_{\xi_7}$ which are not all zero. 

Then necessarily we have $u_1u_4 - u_2u_3 = 0$. 

\bef
{\rm 
For a hyperbolic $(4, 7)$-distribution $(X, D)$, we define 
the {\it characteristic cone} $C \ (\subset D \subset TX)$  by 
\[
\begin{array}{lr}
C := \{ (x, u) \in D \mid 
u = &u_1\xi_1(x) + u_2\xi_2(x) + u_3\xi_3(x)+ u_4\xi_4(x), \quad\quad\quad
\\
&u_1u_4 - u_2u_3 = 0 \}. 
\end{array}
\]
}
\enf

Since the notion of singular curves does not depends on the choice of local frames, we see that 
$C$ does not depend on the choice of local frame $\xi_1, \xi_2, \xi_3, \xi_4$ and is globally well-defined. 
Then we observe the following. 

\bep
\label{necessary-condition}
{\rm (\cite{Montgomery})} \ 
Let $(X, D)$ be a hyperbolic $(4, 7)$-distribution. 
Then there exists uniquely the conformal $(2, 2)$-metric on $D$ with the null-cone field $C \subset TX$ such that, if 
a curve $\gamma : I \to X$ is a $D$-singular curve, then 
$\gamma$ is is tangent to the cone $C$, 
i.e., $\gamma'(t) \in C_{\gamma(t)}$ almost everywhere. 
In particular we have ${\mbox {\it{SVC}}} \subseteq C$. 
\enp

Note that, for each point $x \in X$, 
\[
\begin{array}{l}
C_{x} \setminus \{ 0\}  \cong 
\left\{ 
U = 
\left.
\left(
\begin{array}{cc}
u_1 & u_2 
\\
u_3 & u_4 
\end{array} 
\right) 
\in M_2(\R) 
\right\vert 
\ {\mbox \rm{rank}}(U) = 1 
\right\} 
\\
\ = 
\left.
\left\{ 
\left(
\begin{array}{cc}
\varphi_1
\\
\varphi_2
\end{array} 
\right) 
(\theta_1, \theta_2)
= 
\left(
\begin{array}{cc}
\varphi_1\theta_1 & \varphi_1\theta_2
\\
\varphi_2\theta_1 & \varphi_2\theta_2
\end{array}
\right) 
\right\vert \ 
\left(
\begin{array}{cc}
\varphi_1
\\
\varphi_2
\end{array} 
\right) 
\not= 0, 
(\theta_1, \theta_2) \not= 0
\right\}, 
\end{array}
\]
and that $P(C_{x_0}) \cong \R P^1\times \R P^1$, a torus. 

\section{Prolongations of hyperbolic $(4, 7)$-distributions}
\label{Prolongations of hyperbolic $(4, 7)$-distributions}

To investigate singular curves for a hyperbolic $(4, 7)$-distribution $(X, D)$ further, 
we will use a kind of \lq\lq micro-scope", namely, we will  \lq\lq prolong"
$(X, D)$  along $C$. 

Let $(X, D)$ be a hyperbolic $(4, 7)$-distribution. 
We set 
\[
Z := PC = (C \setminus {\mbox{\rm zero-section}})/\R^\times = \{ (x, [u]) \mid (x, u) \in C, u \not= 0\}/\R^\times. 
\] 
It is regarded as the set of tangential directions in $C$, namely, 
as the set of pairs $(x, \ell)$ of a point $x \in X$ and 
a linear $1$-dimensional subspace $\ell \subset C_x (\subset T_xX)$. 
Then $Z = PC$ is a manifold of dimension $9$ with the natural fibration $\pi : Z \to X$, $\pi(x, [u]) = x$. 
Moreover we define the distribution $E \subset TZ$ 
of rank $3$ by 
\[
E_{(x, [u])} := \{ v \in T_{(x, [u])}Z \mid \pi_*(v) \in \R u \}. 
\]

Let us denote by ${\mathcal E}$ the sheaf of sections to $E$. We define 
${\mathcal E}^{(r+1)} := 
{\mathcal E}^{(r)} + [{\mathcal E}, {\mathcal E}^{(r)}]$ inductively, 
for $r = 1, 2, 3, \dots,$ with ${\mathcal E}^{(1)} := {\mathcal E}$. 
The {\it weak derived system} of $E$ is given by 
${\mathcal E}, {\mathcal E}^{(2)}, {\mathcal E}^{(3)}, \dots$ and 
the {\it small growth} of $E$ at a point $z \in Z$ is defined as the integer vector 
$({\mathrm{rank}}_z {\mathcal E}, \ {\mathrm{rank}}_z {\mathcal E}^{(2)}, \  {\mathrm{rank}}_z {\mathcal E}^{(3)}, 
\dots)$. 

Take a system of local coordinates $x = (x_1, x_2, \dots, x_7)$ around a point $x_0 \in X$. 
Take a point $[u] = [u_1\xi_1(x_0) + u_2\xi_2(x_0) + u_3\xi_1(x_0) + u_4\xi_2(x_0)] \in PC_{x_0}$ 
with $u_1u_4 - u_2u_3 = 0$. Covering the torus $PC_{x_0}$ by 
four coordinate neighbourhoods, we suppose $u_1 \not= 0$. 
(For other cases $u_2 \not= 0, u_3 \not= 0$ or 
$u_4 \not= 0$, similar calculations work well, which we omit them.) 
Then we take the system of local coordinates $(x, \theta, \varphi)$ around $(x_0, [\xi_1(x_0)])$ of $Z$ so that 
${\mathcal E}$ is generated locally in $Z$ by 
\[
\left\{
\begin{array}{l}
\zeta_1 := \dfrac{\partial}{\partial \theta}, 
\\
\zeta_2 := \xi_1 + \theta\xi_2 + \varphi\xi_3 + \theta\varphi\xi_4 + a(x, \theta, \varphi)\dfrac{\partial}{\partial \theta} 
+ b(x, \theta, \varphi)\dfrac{\partial}{\partial \varphi}, 
\\
\zeta_3 := \dfrac{\partial}{\partial \varphi}. 
\end{array}
\right. 
\]
Here $\theta = \theta_2/\theta_1$, $\varphi = \varphi_2/\varphi_1$, where $\varphi_1\theta_1 \not= 0$ 
(see \S \ref{(4,7)-distributions}). 
We will arrange $a, b$ later just to simplify the calculations. 
Note that the growth of $E$ does not depend on the choice of generators, 
in particular, on the choice of $a, b$. 

Now we have $[\zeta_1, \zeta_2] \equiv \xi_2 + \varphi\xi_4 =: \zeta_4, [\zeta_1, \zeta_3] \equiv 0, 
[\zeta_2, \zeta_3] \equiv - \xi_3 - \theta\xi_4 =: \zeta_5, \, {\mbox {\rm mod}}.\, {\mathcal E}$. Up to now we have ${\mathcal E}^{(2)}$, which 
is generated by a distribution $E^{(2)}$ of rank $5$. 

Further we have 
$[\zeta_1, \zeta_4] \equiv 0, [\zeta_1, \zeta_5] \equiv - \xi_4 =: \zeta_6$. 
$[\zeta_2, \zeta_4] \equiv 0$, if we choose $b$ properly, $[\zeta_3, \zeta_4] = \xi_4$, 
$[\zeta_2, \zeta_5] \equiv - \xi_5 - 2\theta\xi_6 - \theta^2\xi_7 - a\xi_4 =:\zeta_7, [\zeta_3, \zeta_5] \equiv 0, \, {\mbox {\rm mod}}. \, {\mathcal E}^{(2)}$. 
Up to now we have ${\mathcal E}^{(3)}$ which is of rank $7$. 
Therefore we see always $E$ has small growth $(3, 5, 7, \dots)$. 

Further we have 
$[\zeta_1, \zeta_6] \equiv 0, [\zeta_1, \zeta_7] \equiv 2(-\xi_6-\theta\xi_7), 
[\zeta_2, \zeta_6] \equiv - \xi_6 - \theta\xi_7 =: \zeta_8, 
[\zeta_3, \zeta_6] \equiv 0, [\zeta_3, \zeta_7] \equiv 0, \, {\mbox {\rm mod}}. \, {\mathcal E}^{(3)}$. 
Moreover we have $[\zeta_2, \zeta_7] \equiv c(x, \theta, \varphi)\xi_7 + d(x, \theta, \varphi)\zeta_8$ 
for some functions $c(x, \theta, \varphi)$ and $d(x, \theta, \varphi)$, ${\mbox {\rm mod}}. 
\, {\mathcal E}^{(3)}$. 
Here it is observed that $c(x, \theta, \varphi)$ is written as 
a polynomial function on $(\theta, \varphi)$. Note that $\zeta_8$ belongs to ${\mathcal E}^{(4)}$ and that 
$[\zeta_2, \zeta_7] \in {\mathcal E}^{(4)}$ if and only if $c(x, \theta, \varphi) = 0$. 
Up to now we have ${\mathcal E}^{(4)}$, which is of rank $8$ if $c(x, \theta, \varphi) = 0$ and of rank $9$ if $c(x, \theta, \varphi) \not= 0$. 

Lastly, by bringing in $[\zeta_1, \zeta_8] = \xi_7$, we have the whole 
$TZ$. We set $\zeta_9 := \xi_7$. 
We remark that, by a proper choice of $a$, we have $[\zeta_2, \zeta_8] \equiv 0, \, 
{\mbox {\rm mod}}. \, {\mathcal E}^{(4)}$. 

Thus we see $E$ has growth $(3, 5, 7, 8, 9)$ 
or $(3, 5, 7, 9)$ depending on whether $c = 0$ or $c \not= 0$. 

When we choose $\xi_2, \xi_3$ or $\xi_4$ instead of $\xi_1$, we have by similar calculations that 
$E$ has small growth $(3, 5, 7, 8, 9)$ 
or $(3, 5, 7, 9)$ for any point $(x_0, [u_0])$ of $Z$. 

\smallskip

Let $(X, D)$ be a hyperbolic $(4, 7)$-distribution, $C \subset TX$ the characteristic cone of $D$, and 
$(Z = PC, E)$ the prolongation of $C$. 
Recall that a direction $[u_0] \in P(C_{x_0})$ at $x_0 \in X$ on $C$ is called {\it of $C_3$-type} 
if $(Z, E)$ has small growth $(3, 5, 7, 8, 9)$ in a neighbourhood of $(x_0, [u_0]) \in Z = PC$ 
(Definition \ref{C_3-type-def}). 

\bef
\label{C_3-and-regular}
{\rm 
A direction $[u_0] \in P(C_{x_0})$ at $x_0 \in X$ on $C$ is called
{\it of regular type}
if $(Z, E)$ has small growth $(3, 5, 7, 9)$ at $(x_0, [u_0])$ in $Z$. 
If $(Z, E)$ has small growth $(3, 5, 7, 8, 9)$ at $(x_0, [u_0])$ in $Z$, but $[u_0]$ is not of $C_3$-type, then 
$[u_0]$ is called {\it of mixed type}. 
}
\enf

Note that, if $[u_0] \in P(C_{x_0})$ is of regular type, then there exists a neighbourhood $U$ 
of $(x_0, [u_0]) \in Z$ in $Z$ such that for any $(x, [u]) \in U$, 
$[u] \in P(C_x)$ is of regular type.  

For each $x_0 \in X$, the locus $\Sigma_{x_0}$ of directions of non-regular type in $P(C_{x_0})$ over $x_0$ 
turns out to be an algebraic set of the torus $P(C_{x_0})$, by the above calculation. 
Therefore we have 

\bel
Let $(X, D)$ be a hyperbolic $(4, 7)$-distribution and $x_0 \in X$. 
Suppose there exists a direction $[u_0] \in P(C_{x_0})$ of $C_3$-type over $x_0$. 
Then $x_0$ is a point of $C_3$-type. Moreover any point in $X$ nearby $x_0$ is of $C_3$-type. 
Thus $C_3$-points of $D$ on $X$ form an open set in $X$. 
\enl

\Proof 
Suppose $[u_0]$ is of $C_3$-type. Then there exists a neighbourhood $U$ of 
$(x_0, [u_0])$ in $Z$ such that $E$ has small growth $(3,5,7,8,9)$ at any point 
$(x, [u]) \in U$. This means that all coefficients of 
the polynomial function $c(x, \theta, \varphi)$ of 
$\theta, \varphi$ vanish in a neighbourhood $V$ of $x_0$ in $X$. Then
any point $x \in V$ is a point of $C_3$-type. 
\QED

\ber
{\rm 
We should remark that a $(4, 7)$-distribution is very degenerate 
at a point of type $C_3$ or at a boundary point of the domain of $C_3$-points 
in the class of all $(4, 7)$-distributions. 
The condition that a $(4, 7)$-distribution is of $C_3$-type is never an \lq\lq open condition\rq\rq. 
}
\enr

\

Let us explain why $C_3$ appears in our story.  

For any point $z \in Z$, define the nilpotent graded Lie algebra 
\[
{\mathfrak g}_{-}  =  {\mathfrak g}_{-5}\oplus{\mathfrak g}_{-4}\oplus{\mathfrak g}_{-3}\oplus
{\mathfrak g}_{-2}\oplus{\mathfrak g}_{-1}, 
\]
where 
\[
\begin{array}{c}
{\mathfrak g}_{-1} := {\mathcal E}_{z}, \ \ {\mathfrak g}_{-2} := ({\mathcal E}^{(2)}/{\mathcal E})_{z}, \ \ 
{\mathfrak g}_{-3} := ({\mathcal E}^{(3)}/{\mathcal E}^{(2)})_{z}, 
\vspace{0.2truecm}
\\
{\mathfrak g}_{-4} := ({\mathcal E}^{(4)}/{\mathcal E}^{(3)})_{z}, \  \ 
{\mathfrak g}_{-5} := (TX/{\mathcal E}^{(4)})_{z}, 
\end{array}
\]
and the Lie algebra multiplication is induced from the 
Lie bracket of sections of $E$. 

Suppose that $[u_0] \in P(C_{x_0})$ 
is a direction of type $C_3$ in the sense of Definition \ref{C_3-and-regular}. 
Then, in a neighbourhood of $z_0 = (x_0, [u_0])$ in $Z$,  the nilpotent graded Lie algebra 
${\mathfrak g}_{-}$ is isomorphic to the nilpotent graded Lie algebra of $C_3$-type 
in the context of Cartan-Tanaka theory(\cite{Tanaka, Yamaguchi}), namely, 
there exists basis $v_1, v_2, v_3$ of ${\mathfrak g}_{-1}$, $v_4, v_5$ of ${\mathfrak g}_{-2}$, 
$v_6, v_7$ of ${\mathfrak g}_{-3}$, $v_8$ of ${\mathfrak g}_{-4}$ and $v_9$ of ${\mathfrak g}_{-5}$ 
such that 
$[v_1, v_2] = v_4, \ [v_1, v_3] = 0, \ [v_2, v_3] = v_5, $
$[v_1, v_4] = 0, \ [v_1, v_5] = v_6, \ [v_2, v_4] = 0, 
[v_2, v_5] = v_7, $
$[v_3, v_4] = - v_6, \ [v_3, v_5] = 0, $
$[v_1, v_6] = 0, \ [v_1, v_7] = 2v_8, \ [v_2, v_6] = v_8, \ [v_2, v_7] = 0, $ 
$[v_3, v_6] = 0, \ [v_3, v_7] = 0, $
$[v_1, v_8] = v_9, \ [v_2, v_8] = 0, \ [v_3, v_8] = 0$. 

Note that ${\mathfrak g}_{-}$ prolongs to the simple Lie algebra of $C_3$-type:
\[
{\mathfrak g} = {\mathfrak g}_{-5}\oplus{\mathfrak g}_{-4}\oplus{\mathfrak g}_{-3}\oplus
{\mathfrak g}_{-2}\oplus{\mathfrak g}_{-1}\oplus{\mathfrak g}_{0}\oplus{\mathfrak g}_{1}\oplus{\mathfrak g}_{2}\oplus
{\mathfrak g}_{3}\oplus{\mathfrak g}_{4}\oplus{\mathfrak g}_{5}, 
\]
which is of dimension $21$, the maximal dimension of symmetries on $(4, 7)$-distributions. 
For details on the relation of Lie algebra of $C_3$-type and $(4, 7)$-distributions, see also \cite{Montgomery} \S 7.12.

\

Now we study on singular curves of the prolongation $(Z, E)$. 
For the system of local coordinates $(x, \theta, \varphi)$, 
take the system of dual coordinates $p, \lambda, \mu$ so that 
$(x, \theta, \varphi; p, \lambda, \mu)$ is a system of local coordinates 
of $T^*Z$. 

The Hamiltonian for singular curves on $(Z, E)$ is given by 
\[
H(x, \theta, \varphi; p. \lambda, \mu) = v_1 H_{\zeta_1} + v_2 H_{\zeta_2} + v_3 H_{\zeta_3}, 
\]
and Hamilton equation on $(x(t), \theta(t), \varphi(t); p(t), \lambda(t), \mu(t))$ is given by 
\[
\dot{x} = \frac{\partial H}{\partial p}, \ 
\dot{\theta} = \frac{\partial H}{\partial \lambda}, 
\ \dot{\varphi} = \frac{\partial H}{\partial \mu}, 
\ \ 
\dot{p} = - \frac{\partial H}{\partial x}, \ \dot{\lambda} = - \frac{\partial H}{\partial \theta}, \ 
\dot{\mu} = - \frac{\partial H}{\partial \varphi}, 
\]
with constraints $H_{\zeta_1} = 0, H_{\zeta_2} = 0, H_{\zeta_3} = 0$ 
and an additional condition $(p(t), \lambda(t), \mu(t)) \not= (0, 0, 0)$. 
Here $v_1, v_2, v_3$ are the control parameters (see \S \ref{Singular curves}). 

\bel
\label{E-singular-curves}
{\rm (1)} If $(Z, E)$ is of $C_3$-type at $z_0 = (x_0, [v]) \in Z$, then 
all $E$-singular paths $\widetilde{\gamma}(t)$ through $z_0$ at $t = 0$ are classified into 
four classes: 
\\
{\rm (a)} $v_2 = 0$ with a singular bi-characteristic $\widetilde{\Gamma}$ in $E^\perp 
\setminus E^{(2)\perp}$,  
\\
{\rm (b)} $v_2 = 0$ or $v_3 = 0$ with a singular bi-characteristic $\widetilde{\Gamma}$ in $E^{(2)\perp} \setminus E^{(3)\perp}$, 
\\
{\rm (c)} $(v_1, v_2) = (0, 0)$ with a singular bi-characteristic $\widetilde{\Gamma}$ in 
$E^{(3)\perp} \setminus E^{(4)\perp}$, 
\\
{\rm (d)} $v_1 = 0$ with a singular bi-characteristic $\widetilde{\Gamma}$ in $E^{(4)\perp} 
\setminus {\mbox{\rm (zero-section)}}$. 

{\rm (2)} If $(Z, E)$ is of regular type at $z_0 = (x_0, [v]) \in Z$ 
where the growth of $E$ is $(3, 5, 7, 9)$, 
all $E$-singular paths $\Gamma(t)$ through $z_0$ at $t = 0$ are classified into 
three classes: 
\\
{\rm (a)} $v_2 = 0$ with a singular bi-characteristic $\widetilde{\Gamma}$ in $E^\perp\setminus E^{(2)\perp}$,  
\\
{\rm (b)} $v_2 = 0$ or $v_3 = 0$ with a singular bi-characteristic $\widetilde{\Gamma}$ in $E^{(2)\perp}\setminus E^{(3)\perp}$, 
\\
{\rm (c)} $(v_1, v_2) = (0, 0)$ with a singular bi-characteristic $\widetilde{\Gamma}$ in 
$E^{(3)\perp}\setminus E^{(4)\perp}$. 
\enl

If $(Z, E)$ is of $C_3$-type (resp. of regular type) at $(x_0, [v])$, 
then the singular velocity cone of $(Z, E)$ at $(x_0, [v])$ is depicted as below left (resp. below right): 
\begin{center}
\includegraphics[width=4truecm, height=3truecm, clip, 
]{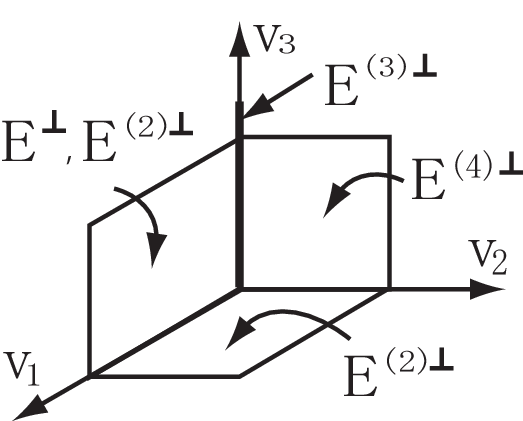} 
\qquad 
\includegraphics[width=4truecm, height=3truecm, clip, 
]{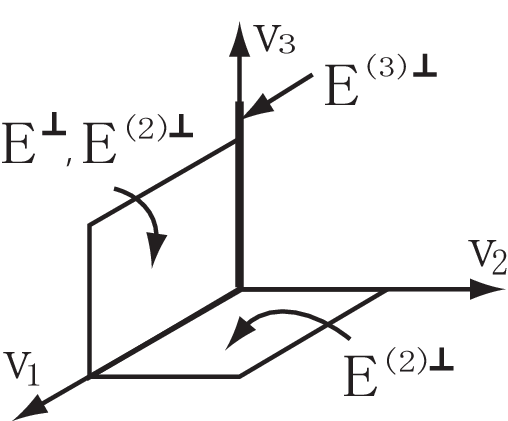} 
\\
$C_3$ \hspace{4truecm} regular
\end{center}

\

\noindent
{\it Proof of Lemma \ref{E-singular-curves}:} \ 
Let 
$\widetilde{\gamma}(t) = (x(t), \theta(t), \varphi(t))$ be 
a non-constant $E$-singular curve and 
$\widetilde{\Gamma}(t) = (x(t), \theta(t), \varphi(t); p(t), \lambda(t), \mu(t))$ 
an associated singular bi-characteristic curve in $T^*Z$ to $\widetilde{\gamma}(t)$. 

From the bracket relations of $E$, we have, modulo $\langle H_{\zeta_1}, H_{\zeta_2}, H_{\zeta_3}\rangle$, 
\begin{align*}
\frac{d}{dt}H_{\zeta_1} & =  
v_1H_{[\zeta_1, \zeta_1]} + v_2H_{[\zeta_2, \zeta_1]} + v_3H_{[\zeta_3, \zeta_1]} 
\equiv - v_2H_{\zeta_4}. 
\\
 \frac{d}{dt}H_{\zeta_2} & =  
 v_1H_{[\zeta_1, \zeta_2]} + v_2H_{[\zeta_2, \zeta_2]} + v_3H_{[\zeta_3, \zeta_2]} 
\equiv v_1H_{\zeta_4} - v_3H_{\zeta_5}, 
\\
\frac{d}{dt}H_{\zeta_3} & =  v_1H_{[\zeta_1, \zeta_3]} + v_2H_{[\zeta_2, \zeta_3]} + v_3H_{[\zeta_3, \zeta_3]} 
\equiv v_2H_{\zeta_5}. 
\end{align*}
Here  $\langle H_{\zeta_1}, H_{\zeta_2}, H_{\zeta_3}\rangle$ means the ideal generated by 
$H_{\zeta_1}, H_{\zeta_2}, H_{\zeta_3}$, composed with the bi-characteristic $\widetilde{\Gamma} 
: (\R, 0) \to T^*Z$, over the ring of $C^\infty$ functions on $t$. 

Moreover we have, modulo $\langle H_{\zeta_1}, H_{\zeta_2}, H_{\zeta_3}, H_{\zeta_4}, H_{\zeta_5}\rangle$, 
\begin{align*}
\frac{d}{dt}H_{\zeta_4} & =  v_1H_{[\zeta_1, \zeta_4]} + v_2H_{[\zeta_2, \zeta_4]} + v_3H_{[\zeta_3, \zeta_4]} 
\equiv - v_3H_{\zeta_6}, 
\\
\frac{d}{dt}H_{\zeta_5} & =  v_1H_{[\zeta_1, \zeta_5]} + v_2H_{[\zeta_2, \zeta_5]} + v_3H_{[\zeta_3, \zeta_5]} 
\equiv v_1H_{\zeta_6} + v_2H_{\zeta_7}. 
\end{align*}

Further we have, modulo $\langle H_{\zeta_1}, H_{\zeta_2}, H_{\zeta_3}, H_{\zeta_4}, H_{\zeta_5}, H_{\zeta_6}, H_{\zeta_7}\rangle$, 
\begin{align*}
\frac{d}{dt}H_{\zeta_6} & =  v_1H_{[\zeta_1, \zeta_6]} + v_2H_{[\zeta_2, \zeta_6]} + v_3H_{[\zeta_3, \zeta_6]} 
\equiv v_2H_{\zeta_8}, 
\\
\frac{d}{dt}H_{\zeta_7} & =  v_1H_{[\zeta_1, \zeta_7]} + v_2H_{[\zeta_2, \zeta_7]} + v_3H_{[\zeta_3, \zeta_7]} 
\equiv 2v_1H_{\zeta_8} + v_2cH_{\zeta_9}. 
\end{align*}

Furthermore we have, modulo 
$\langle H_{\zeta_1}, H_{\zeta_2}, H_{\zeta_3}, H_{\zeta_4}, H_{\zeta_5}, H_{\zeta_6}, H_{\zeta_7}, H_{\zeta_8}\rangle$, 
\begin{align*}
\frac{d}{dt}H_{\zeta_8} & =  v_1H_{[\zeta_1, \zeta_8]} + v_2H_{[\zeta_2, \zeta_8]} + v_3H_{[\zeta_3, \zeta_8]} 
\equiv v_1H_{\zeta_9}. 
\end{align*}
Here recall that $[\zeta_2, \zeta_7] \equiv c\zeta_9, \, {\mbox {\rm mod}}.\, {\mathcal E}^{(3)}$ 
for a function $c$ on $(x, \theta, \varphi)$, 
and that, 
for the case (1), namely, for a point $(x, \theta, \varphi)$ where the growth of $E$ is $(3, 5, 7, 8, 9)$ locally, we have $c(x, \theta, \varphi) = 0$ (see \S \ref{Prolongations of hyperbolic $(4, 7)$-distributions}). 

Now, from the constraints $H_{\zeta_1} = 0, H_{\zeta_2} = 0$ and $H_{\zeta_3} = 0$, necessarily we have along $\widetilde{\Gamma}$, 
$v_2H_{\zeta_4} = 0, v_1H_{\zeta_4} - v_3H_{\zeta_5} = 0$ and $v_2H_{\zeta_5} = 0$. 

For any parameter $t$ where $(H_{\zeta_4}, H_{\zeta_5}) \not= (0, 0)$, 
we have $v_2(t) = 0$. 
Then the ratio $v_1(t) : v_3(t)$ is determined as $H_{\zeta_5} : H_{\zeta_4}$. Then 
any solution of the ordinary differential equation 
of the vector field $H_{\zeta_5}\overrightarrow{H_{\zeta_1}} + H_{\zeta_4}\overrightarrow{H_{\zeta_3}}$ over $T^*X$ satisfying 
the initial conditions $H_{\zeta_1} = 0, H_{\zeta_2}  = 0, H_{\zeta_3} = 0$ gives 
an $E$-singular curve in $E^\perp \setminus E^{(2)\perp}$ near $t = 0$. 
This corresponds to the case (1)(a) or (2)(a). 

Next suppose $(H_{\zeta_4}, H_{\zeta_5}) = (0, 0)$ on an interval $J$. 
Then $v_3H_{\zeta_6} = 0, v_1H_{\zeta_6} + v_2H_{\zeta_7} = 0$ on $J$. 

For any parameter $t$ where $H_{\zeta_6} \not= 0$, 
we have $v_3(t) = 0$. Then the ration 
$v_1(t) : v_2(t)$ is determined as $-H_{\zeta_7} : H_{\zeta_6}$. 
Consider the vector field 
$\widetilde{\zeta} := -H_{\zeta_7}\overrightarrow{H_{\zeta_1}} + H_{\zeta_6}\overrightarrow{H_{\zeta_2}}$ 
over $T^*Z$. 
Then we have that each $\widetilde{\zeta}H_{\zeta_i}$ belongs to the ideal generated by 
$H_{\zeta_1}, \dots, H_{\zeta_5}$ over $\widetilde{\Gamma}(t)$, for $i = 1, 2, 3, 4, 5$. 
Then $\widetilde{\Gamma}(t)$ is obtained, up to parametrisation, 
as a solution of the ordinary differential equation of the vector field $\widetilde{\zeta}$ 
with an initial conditions 
$H_{\zeta_1} = H_{\zeta_2} = H_{\zeta_3} = H_{\zeta_4} = H_{\zeta_5} = 0, H_{\zeta_6} \not= 0$. 
Therefore $\widetilde{\Gamma}(t)$ is contained in $E^{(2)\perp} \setminus E^\perp$. 

For any parameter $t$ where $H_{\zeta_6} = 0, H_{\zeta_7} \not= 0$, 
we have $v_2(t) = 0$. For any functions $A(x, \theta, \varphi), B(x, \theta, \varphi)$ 
we consider a vector field $\widetilde{\zeta} := A\overrightarrow{H_{\zeta_1}} + B\overrightarrow{H_{\zeta_3}}$ over $T^*Z$. Then 
we have $\widetilde{\zeta}H_{\zeta_i}$ belongs to the ideal generated by 
$H_{\zeta_1}, \dots, H_{\zeta_6}$ over $\widetilde{\Gamma}$, for $i = 1, 2, 3, 4, 5, 6$. 
Then $\widetilde{\Gamma}(t)$ is obtained, up to parametrisation, 
as a solution of the ordinary differential equation of the vector field $\widetilde{\zeta}$ 
with an initial conditions 
$H_{\zeta_1} = H_{\zeta_2} = H_{\zeta_3} = H_{\zeta_4} = H_{\zeta_5} = 0, H_{\zeta_6} = 0, H_{\zeta_7} \not= 0$. 
Therefore $\widetilde{\Gamma}(t)$ is contained in $E^{(2)\perp} \setminus E^\perp$. 

These cases correspond to the case (1)(b) or (2)(b). 

Next suppose $(H_{\zeta_6}, H_{\zeta_7}) = (0, 0)$ on an interval $K$. 
Then $v_2H_{\zeta_8} = 0, 2v_1H_{\zeta_8} + v_2cH_{\zeta_9} = 0$ on $K$. 
If $H_{\zeta_8} \not= 0$, then we have $v_1(t) = v_2(t) = 0$. 
Set $\widetilde{\zeta} := \overrightarrow{H_{\zeta_3}}$. 
Then we have $\widetilde{\zeta}H_{\zeta_i}$ belongs to the ideal generated by 
$H_{\zeta_1}, \dots, H_{\zeta_7}$ over $\widetilde{\Gamma}$, for $i = 1, 2, 3, 4, 5, 6, 7$. 
Then $\widetilde{\Gamma}(t)$ is obtained, up to parametrisation, 
as a solution of the ordinary differential equation of the vector field $\widetilde{\zeta}$ 
with an initial conditions 
$H_{\zeta_1} = H_{\zeta_2} = H_{\zeta_3} = H_{\zeta_4} = H_{\zeta_5} = 0, H_{\zeta_6} = 0, H_{\zeta_7} = 0, 
H_{\zeta_8} \not= 0$, which 
is an $E$-singular bi-characteristic in $E^{(3)\perp} \setminus E^{(2)\perp}$. 
This corresponds to the case (1)(c) or (2)(c). 

Suppose $H_{\zeta_8} = 0$ on an interval $L$. 
Then we have $v_1H_{\zeta_9} = 0$. Since 
$(p(t), \lambda(t), \mu(t)) \not= (0, 0, 0)$ for corresponding bi-characteristic 
$\widetilde{\Gamma}(t)$, we have $H_{\zeta_9}(t) \not= 0$. 
For any functions $A(x, \theta, \varphi), B(x, \theta, \varphi)$ 
we consider a vector field $\widetilde{\zeta} := AH_{\zeta_2} + BH_{\zeta_3}$ over $T^*Z$. 
Then we have $\widetilde{\zeta}H_{\zeta_i}$ belongs to the ideal generated by 
$H_{\zeta_1}, \dots, H_{\zeta_8}$ over $\widetilde{\Gamma}$, for $i = 1, 2, 3, 4, 5, 6, 7, 8$. 
Then $\widetilde{\Gamma}(t)$ is obtained, up to parametrisation, 
as a solution of the ordinary differential equation of the vector field $\widetilde{\zeta}$ 
with an initial conditions 
$H_{\zeta_1} = H_{\zeta_2} = H_{\zeta_3} = H_{\zeta_4} = H_{\zeta_5} = 0, H_{\zeta_6} = 0, H_{\zeta_7} = 0, 
H_{\zeta_8} = 0, H_{\zeta_9} \not= 0$, which 
is an $E$-singular bi-characteristic in $E^{(4)\perp} \setminus {\mbox{\rm (zero-section)}}$. 
This corresponds to the case (1)(d). 

For a point $(x, \theta, \varphi)$ where the growth of $E$ is $(3, 5, 7, 9)$, 
we have $c(x, \theta, \varphi) \not= 0$. Then, if $v_2(t) \not= 0$, 
then we have $(H_{\zeta_8}, H_{\zeta_9}) = (0, 0)$ and 
then $(p(t), \lambda(t), \mu(t)) = (0, 0, 0)$. 
This means that there does not exist any $E$-singular curve around that point. Therefore $E$-singular paths of type (d) do not appear 
for the case (2). 
\QED

\section{Proof of main theorems}
\label{Proof of Theorem}

Let $(X, D)$ be a hyperbolic $(4, 7)$-distribution, and $C$ the characteristic cone field of $D$. 

\bel
\label{D-singular-and-E-singular}
Let $\gamma : (\R, 0) \to X$ be a germ of $C^\infty$ $C$-integral curve with $\gamma'(0) \not= 0$. 
Then we have 

{\rm (1)} $\gamma$  lifts uniquely to an $E$-integral curve 
$\widetilde{\gamma} : (\R, 0) \to Z$. 

{\rm (2)} 
$\gamma$ is a $D$-singular curve if and only if 
the lift $\widetilde{\gamma}$ of $\gamma$ is an $E$-singular curve with a bi-characteristic 
$\widetilde{\Gamma}$ in $(E^{(3)})^\perp \setminus \mbox{\rm (zero-section)}$, i.e., 
$\widetilde{\gamma}$ is an $E$-integral curve and an $E^{(3)}$-singular curve. 
\enl

\Proof
(1) The $C$-integral curve-germ $\gamma$ induces the curve-germ $\gamma' : (\R, 0) \to TC \setminus ({\mbox{\rm zero-section}})$ 
and the curve-germ $\widetilde{\gamma} : (\R, 0) \to PC = Z$ by $\widetilde{\gamma}(t) := (\gamma(t), [\gamma'(t)])$. 
Then $\pi_*(\widetilde{\gamma}'(t)) = \gamma'(t) \in C_{\gamma(t)}$. Here $\pi : Z = PC \to X$ is the canonical projection. 
Therefore $\widetilde{\gamma}'(t) \in E_{\widetilde{\gamma}(t)}$. 
Hence $\widetilde{\gamma}$ is an $E$-integral curve. 

(2) Recalling the situation, we write down Hamiltonian equation for $(X, D)$ on the system of local coordinates 
$(x; p)$ of $T^*X$: 
(*) $\dot{x} = \frac{\partial H_X}{\partial p}, \  \ \dot{p} = -\frac{\partial H_X}{\partial x}$, 
with constraints $H_{\xi_1} = 0,  H_{\xi_2} = 0,  H_{\xi_3} = 0,  H_{\xi_4} = 0,  p(t) \not= 0$, where 
$H_X = u_1 H_{\xi_1} + u_2 H_{\xi_2} + u_3 H_{\xi_3} + u_4 H_{\xi_4}$. See \S \ref{(4,7)-distributions}. 

Also we write down Hamiltonian equation for $(Z, E)$ on the system of local coordinates 
$(x, \theta, \varphi; p, \lambda, \mu)$ of $T^*Z$: 
(**) $\dot{x} = \frac{\partial H_Z}{\partial p}, \ \dot{\theta} = \frac{\partial H_Z}{\partial \lambda}, \ \dot{\varphi} = \frac{\partial H_Z}{\partial \mu}; \ 
\dot{p} = -\frac{\partial H_Z}{\partial x}, \ \dot{\lambda} = -\frac{\partial H_Z}{\partial \theta}, \ \dot{\mu} = -\frac{\partial H_Z}{\partial \varphi}$, 
with constraints $H_{\zeta_1} = 0, H_{\zeta_2} = 0, H_{\zeta_3} = 0, (p(t), \lambda(t), \mu(t)) \not= (0, 0, 0)$, 
where $H_Z = v_1H_{\zeta_1} + v_2 H_{\zeta_2} + v_3 H_{\zeta_3}, \ 
H_{\zeta_1} = \lambda, \ H_{\zeta_2} = H_{\xi_1} + \theta  H_{\xi_2} + \varphi  H_{\xi_3} + \theta\varphi  H_{\xi_4} + a\lambda + b\mu, 
\ H_{\zeta_3} = \mu$. See \S \ref{Prolongations of hyperbolic $(4, 7)$-distributions}. 

Let $\gamma(t) = x(t)$ be a $D$-singular curve with bi-characteristic $\Gamma(t) = (x(t); p(t))$ satisfying (*). 
Then $\dot{x} = u_1\xi_1 + u_2\xi_2 + u_3\xi_3 + u_4\xi_4$. 
Since $\gamma(t)$ is not constant $u_1, u_2, u_3, u_4$ are not all zero. 
Now suppose $u_1 \not= 0$. 
For other cases, it is proved the required result by a similar argument. 
Now we have $\theta = u_2/u_1, \varphi = u_3/u_1$. Since 
$u_1u_4 - u_2u_3 = 0$, we have $u_4/u_1 = \theta\varphi$. Set $v_1 = \dot{\theta}, v_2 = u_1, v_3 = \dot{\varphi}, 
\lambda(t) = 0, \mu(t) = 0$. 
By the equality $u_1H_{\xi_5} + u_2 H_{\xi_6} = 0, u_1 H_{\xi_6} + u_2 H_{\xi_7} = 0$, 
we have $H_{\xi_5} + 2\theta H_{\xi_6} + \theta^2 H_{\xi_7} = 0$ (See \S \ref{(4,7)-distributions}). 
Therefore, by the arguments in \S \ref{Prolongations of hyperbolic $(4, 7)$-distributions}, 
we see $\widetilde{\gamma}(t) = (x(t), \theta(t), \varphi(t))$ is an $E$-singular curve
with bi-characteristic $\widetilde{\Gamma}(t) = (x(t), \theta(t), \varphi(t); p(t), 0, 0)$ in $E^{(3)\perp}$ 
 satisfying (**). 

Conversely, let $\widetilde{\gamma}(t) = (x(t), \theta(t), \varphi(t))$ be an $E$-singular curve with a bi-characteristic 
$\widetilde{\Gamma}(t) = (x(t), \theta(t), \varphi(t); p(t), \lambda(t), \mu(t))$ in $(E^{(3)})^\perp \setminus \{ 0\}$. 
Since $\widetilde{\Gamma}(t) \in (E^{(3)})^\perp \setminus \{ 0\}$, we have $\lambda = 0, \mu = 0, 
H_{\xi_1} = 0, H_{\xi_2} = 0, H_{\xi_3} = 0, H_{\xi_4} = 0$. 
In particular we have $p(t) \not= 0$. By $\dot{x} = \frac{\partial H_Z}{\partial p}, \dot{p} = -\frac{\partial H_Z}{\partial x}$ 
we have $\dot{x} = \frac{\partial H_X}{\partial p}, \dot{p} = -\frac{\partial H_X}{\partial x}$, by setting 
$u_1 = v_2, u_2 = v_2\theta, u_3 = v_2\varphi, u_4 = v_3\theta\varphi$, and we see $\gamma(t) = x(t)$ 
is a $D$-singular curve with the bi-characteristic $\Gamma(t) = (x(t); p(t))$. 
Similar argument works well for other coordinate neighbourhoods of $Z$, which we omit the details. 
\QED

\
 
Thus, to study singular curves of $(X, D)$, 
we are interested in the singular velocity cone for $(Z, E)$ at a point $z_0 = (x_0, [u_0])$, 
with information on associated bi-characteristics, i.e., 
the set of $(v_1, v_2, v_3)$ such that 
there exists an $E$-integral curve $\widetilde{\gamma} : (\R, 0) \to (Z, z_0)$ 
with the initial vector $\widetilde{\gamma}'(0) = v_1\zeta_1(z_0) + v_2\zeta_2(z_0) + v_3\zeta_3(z_0)$, with information whether it is $E^{(i)}$-singular as well, $i = 1, 2, 3, 4$. 
Then, using  Lemma \ref{E-singular-curves}, we show

\bet 
\label{main-theorem3}
{\rm (1)} If $[u_0] \in P(C_{x_0})$ is of $C_3$-type, then 
there exists a $D$-singular path $\gamma : (\R, 0) \to X$ with $\gamma(0) = x_0, \gamma'(0) = u_0$ 
and therefore $u_0$ belongs to the singular velocity cone ${\mbox {\it{SVC}}}$: $u_0 \in {\mbox {\it{SVC}}}$. 

{\rm (2)} If $[u_0] \in P(C_{x_0})$ is of regular type, then 
there does not exist any $D$-singular path $\gamma : (\R, 0) \to X$ with $\gamma(0) = x_0, \gamma'(0) = u_0$ 
and therefore $u_0 \in C$ never belong to the singular velocity cone ${\mbox {\it{SVC}}}$: $u_0 \not\in {\mbox {\it{SVC}}}$. 
\ent

\Proof
(1) Let $[u_0] \in P(C_{x_0})$ be a direction of $C_3$-type. 
Then by Lemma \ref{E-singular-curves}(1), 
there exists an $E$-singular curve $\widetilde{\gamma} : (\R, 0) \to Z$ with a 
bi-characteristic $\widetilde{\Gamma}$ in $E^{(3)\perp} \setminus {\mbox{\rm (zero-section)}}$ 
such that $\widetilde{\gamma}(0) = (x_0, [u_0])$. 
Then $\gamma = \pi\circ \widetilde{\gamma} : (\R, 0) \to 
X$ is a $D$-singular path with $[\gamma'(0)] = [u_0]$ by Lemma \ref{D-singular-and-E-singular}(2). 
Therefore we see that $u_0 \in {\mbox {\it{SVC}}}_{x_0}$. 

(2) Let $[u_0] \in P(C_{x_0})$ be a direction of regular-type. 
Then by Lemma \ref{E-singular-curves}(2), 
there does not exist any $E$-singular curve $\widetilde{\gamma} : (\R, 0) \to Z$ with a 
bi-characteristic $\widetilde{\Gamma}$ in $E^{(3)\perp}$ such that $\widetilde{\gamma}(0) = (x_0, [u_0])$. 
Therefore there does not exist any $D$-singular path 
$\gamma : (\R, 0) \to X$ with $\gamma(0) = x_0, \gamma'(0) = u_0$. 
Then by Lemma \ref{D-singular-and-E-singular}(2), we have $u_0 \not\in {\mbox {\it{SVC}}}$. 
\QED

\

\noindent
{\it Proof of Theorem \ref{main-theorem}:} 
(1) Let $(X, D)$ be a hyperbolic $(4, 7)$-distribution and $[u] \in P(C_x)$ a direction of type $C_3$. 
By Theorem \ref{main-theorem3} (1), there exists uniquely a $D$-singular curve $\gamma : (\R, 0) \to X$ with $\gamma(0) = x, \gamma'(0) = u$ up to parametrisation. 

(2) 
Let $(X, D)$ be a hyperbolic $(4, 7)$-distribution of type $C_3$. Then, for any $x \in X$, 
for any $u \in C_{x_0} \setminus \{ 0\}$, 
the direction $[u] \in P(C_{x})$ is of type $C_3$. 
By Theorem \ref{main-theorem3} (1), $u \in {\mbox {\it{SVC}}}_{x}$. Note that 
$0 \in {\mbox {\it{SVC}}}$, since constant curves are $D$-singular. Therefore we have ${\mbox {\it{SVC}}} = C$, 
which completes the proof. 
\QED

\

\noindent
{\it Proof of Theorem \ref{main-theorem2}:} 
Let $F$ be the component of the singular velocity cone of $E$ with 
bi-characteristic in $E^{(3)\perp}$ (see Lemma \ref{E-singular-curves}). 

(1) Let $c : I \to Z$ be a $F$-integral curve with $\pi\circ c(t) \not= 0$. Set 
$c'(t) = v_2(t)\zeta_2(c(t)) + v_3(t)\zeta_3(c(t))$. Since $(\pi\circ c)'(t) \not= 0$, we have 
$v_2(t) \not= 0$. Then, locally on $Z$, there exist $C^\infty$ functions $A(x, \theta, \varphi)$, 
$B(x, \theta, \varphi)$ such that $v_2(t) = A(c(t)), v_3(t) = B(c(t))$. 
Then a singular bi-characteristic associated to $c(t)$ is obtained as a solution curve 
of ordinary differential equation defined by the vector field $A\overrightarrow{H_{\zeta_2}} + 
B\overrightarrow{H_{\zeta_3}}$ as in an argument in the proof of Lemma \ref{E-singular-curves}. 
Thus we see that $c(t)$ is an $E$-singular path belonging to the class (1)(d) of Lemma \ref{E-singular-curves}. 
Then, by Lemma \ref{D-singular-and-E-singular}(2), we see $\pi\circ c$ is a $D$-singular path. 

(2) 
Let $(x, [u]) \in Z$ and $\gamma : (\R, 0) \to X$ be a $D$-singular path with $\gamma(0) = x$ and 
$[\gamma'(t)] = [u]$. Take the unique lift $\widetilde{\gamma} : (\R, 0) \to Z$ of $\gamma$. 
Then by Lemma \ref{D-singular-and-E-singular}(2), $\widetilde{\gamma}$ is an $E$-singular 
path which belongs to the class (1)(d) of Lemma \ref{D-singular-and-E-singular}(2). 
Thus we see $\widetilde{\gamma}$ is a $F$-integral curve. 
\QED

\section{Isotropic Grassmannian on $6$-dimensional symplectic space}
\label{Isotropic-Grassmannian}

Here we give examples of hyperbolic $(4, 7)$-distributions. 

Let $\Omega$ be the symplectic form on the vector space $\R^6$ with basis 
$e_1, \dots, e_6$ having the matrix representation 
\[
\Omega = 
\left(
\begin{array}{cccccc}
0 & 0 & 0 & 0 & 0 & 1
\\
0 & 0 & 0 & 0 & 1 & 0
\\
0 & 0 & 0 & 1 & 0 & 0
\\
0 & 0 & -1 & 0 & 0 & 0
\\
0 & -1 & 0 & 0 & 0 & 0
\\
-1 & 0 & 0 & 0 & 0 & 0
\end{array}
\right)
\]
The symplectic group ${\textrm{Sp}}(6, \R)$ acts transitively on the isotropic Grassmannian ${\mathcal F}_2$ of 
$2$-dimensional isotropic subspaces in $\R^6$. We set $X = {\mathcal F}_2$. 
Any $V_2 \in X$ near $\langle e_1, e_2\rangle_{\R}$ has a basis $h_1, h_2$ of form 
\[
\left\{
\begin{array}{rcccl}
h_1 & = & e_1 &     & + x_1^3 r_3 + x_1^4 e_4 + x_1^5 e_5 + x_1^6 e_6,  
\\
h_2 & = &       &e_2 & + x_2^3r_3 + x_2^4e_4 + x_2^5e_5 + x_2^6e_6,  
\end{array}
\right.
\]
with $\Omega(h_1, h_2) = x_2^6 + x_1^3x_2^4 - x_1^4x_2^3 - x_1^5 = 0$. 
Thus the homogeneous space $X$ is of dimension $7$ 
and has a system of local coordinates $x_1^3, x_1^4, x_1^5, x_1^6, x_2^3, x_2^4, x_2^5$
in a neighbourhood $U$ of $\langle e_1, e_2\rangle_{\R}$. 
The canonical distribution $D \subset TX$ is defined such that, for any $(x, v) \in TX$ with
$x \in U$, $v \in D_{x, v}$ if and only if there exists a $C^\infty$ curve $V_2(t)$ on $X$ starting from 
$x$ at $t = 0$ and $V_2'(t)$ is contained in the skew-orthogonal $V_2(t)^\perp$ of $V_2(t)$ for 
$\Omega$. This means that 
there exists a $C^\infty$ family $(h_1(t), h_2(t))$ with $x = 
\langle h_1(0), h_2(0)\rangle_\R$ and with 
\[
\left\{
\begin{array}{rcl}
\Omega(h_1(t), h_1'(t)) & = & {x_1^6}' + x_1^3{x_1^4}' - x_1^4{x_1^3}' = 0,  
\\
\Omega(h_1(t), h_2'(t)) & = & {x_2^6}' + x_1^3{x_2^4}' - x_1^4{x_2^3}' = 0, 
\\
\Omega(h_2(t), h_1'(t)) & = & {x_1^5}' + x_2^3{x_1^4}' - x_2^4{x_1^3}' = 0,  
\\
\Omega(h_2(t), h_2'(t)) & = & {x_2^5}' + x_2^3{x_2^4}' - x_2^4{x_2^3}' = 0. 
\end{array}
\right.
\]
Note that $x_2^6 = x_1^5 -  x_1^3x_2^4 + x_1^4x_2^3$. Then it is proved that 
the above second equality follows from other three equalities, 
that $D$ is locally defined by a Pfaff system
\[
dx_1^6 + x_1^3dx_1^4 - x_1^4dx_1^3 = 0, \ \ 
dx_1^5 + x_2^3dx_1^4 - x_2^4dx_1^3 = 0, \ \ 
dx_2^5 + x_2^3dx_2^4 - x_2^4dx_2^3 = 0, 
\]
$D$ is of rank $4$ and that $D$ has a local frame
\[
\begin{array}{c}
\xi_1 = \frac{\partial}{\partial x_1^3} + x_2^4\frac{\partial}{\partial x_1^5} + x_1^4\frac{\partial}{\partial x_1^6}, \ 
\xi_2 = \frac{\partial}{\partial x_2^3} + x_2^4\frac{\partial}{\partial x_2^5}, 
\vspace{0.2truecm}
\\
\xi_3 = \frac{\partial}{\partial x_1^4} - x_2^3\frac{\partial}{\partial x_1^5} - x_1^3\frac{\partial}{\partial x_1^6}, \ 
\xi_4 = \frac{\partial}{\partial x_2^4} - x_2^3\frac{\partial}{\partial x_2^5}
\end{array}
\]
over $U$. Then we have $[\xi_1, \xi_2] = 0, \ [\xi_1, \xi_3] = -2\frac{\partial}{\partial x_1^6} =: \xi_5, \ 
[\xi_1, \xi_4] = - \frac{\partial}{\partial x_1^5} =: \xi_6, \ [\xi_2, \xi_3] = - \frac{\partial}{\partial x_1^5}, \ 
[\xi_2, \xi_4] = - 2\frac{\partial}{\partial x_2^5} =: \xi_7, \ [\xi_3, \xi_4] = 0$, and that 
$D$ is a hyperbolic $(4, 7)$-distribution. 
Moreover the prolongation $(Z, E)$ of the characteristic cone $C \subset D$ is 
isomorphic to the isotropic flag manifold ${\mathcal F}_{1,2,3}$ on $(\R^6, \Omega)$ with canonical distribution $E$, 
which has growth $(3, 5, 7, 8, 9)$. 

\ber
\label{Engel-foliation} 
{\rm 
Suppose a $(4, 7)$-distribution $D$ has a frame $\xi_1, \xi_2, \xi_3, \xi_4$ such that 
$[\xi_1, \xi_2] = 0, [\xi_3, \xi_4] = 0, [\xi_1, \xi_4] = [\xi_2, \xi_3], [\xi_i, [\xi_j, \xi_k]] = 0\, (1 \leq i, j, k \leq 4)$, 
as in the above model of hyperbolic $(4, 7)$-distribution $(X, D)$. 
Then the distribution $F \subset D$ generated by $\zeta_2$ and $\zeta_3$ of rank $2$ 
of the prolongation $(Z, E)$ (see Theorem \ref{main-theorem2} and \S \ref{Prolongations of hyperbolic $(4, 7)$-distributions}) 
have the small derived systems $F^{(2)} = \langle \zeta_2, \zeta_3, \zeta_5\rangle$ and $F^{(3)} = \langle \zeta_2, \zeta_3, \zeta_5, \zeta_7\rangle$ which turns to be complete integrable. 
Then $F^{(3)}$ induces an Engel foliation on $Z$ i.e., a foliation consisting of four dimensional leaves $L$ with Engel distribution $F\vert_L$. 
}
\enr

\

Based on the above model, we give an example of families of hyperbolic $(4, 7)$-distributions. 

Let $\R^7 = \{ (x, y, s, t, z, w, u) \}$, and consider the concrete distribution $D$ on $X = \R^7$, 
with parameter $\varepsilon \in \R$, 
generated by 
\[
\begin{array}{c}
\xi_1 =  \frac{\partial}{\partial x} + (w + \varepsilon t)\frac{\partial}{\partial s} + y\frac{\partial}{\partial t}, \quad
\xi_2  =  \frac{\partial}{\partial z} + w\frac{\partial}{\partial u}, 
\vspace{0.2truecm}
\\
\xi_3  =  \frac{\partial}{\partial y} - (z - \varepsilon t)\frac{\partial}{\partial s} - x\frac{\partial}{\partial t}, \quad
\xi_4  =  \frac{\partial}{\partial w} - z\frac{\partial}{\partial u}. 
\end{array}
\]
We have 
\[
\begin{array}{c}
[\xi_1, \xi_2] = 0, \ \ [\xi_1, \xi_3] = - 2\frac{\partial}{\partial t} + \varepsilon(x+y)\frac{\partial}{\partial s} =: \xi_5, 
\ \ 
[\xi_1, \xi_4] = - \frac{\partial}{\partial s} =: \xi_6, 
\\
{[\xi_2, \xi_3]} = - \frac{\partial}{\partial s}, \ \ 
[\xi_2, \xi_4] = -2\frac{\partial}{\partial u} =: \xi_7, \ \  [\xi_3, \xi_4] = 0. 
\end{array}
\]
Then $(X, D)$ is a hyperbolic $(4, 7)$-distribution and $(\xi_1, \xi_2, \xi_3, \xi_4)$ is an adapted 
frame. Moreover we have 
\[
\begin{array}{l}
[\xi_1, \xi_5] = -3\varepsilon \xi_6, \  [\xi_1, \xi_6] = 0, \  [\xi_1, \xi_7] = 0, \  
[\xi_2, \xi_5] = 0, \  [\xi_2, \xi_6] = 0, 
\\
{[\xi_2, \xi_7]} = 0, \ 
[\xi_3, \xi_5] = - 3\varepsilon \xi_6, \ [\xi_3, \xi_6] = 0, \ [\xi_3, \xi_7] = 0, \ 
[\xi_4, \xi_5] = 0, \ 
\\
{[\xi_4, \xi_6]} = 0, \ [\xi_4, \xi_7] = 0. 
\end{array}
\]

The frame of $E$ as in \S \ref{Prolongations of hyperbolic $(4, 7)$-distributions} 
is given by 
\[
\zeta_1 = \frac{\partial}{\partial \theta}, \quad \zeta_2 = \xi_1 + \theta\xi_2 + \varphi\xi_3 + \theta\varphi\xi_4, 
\quad 
\zeta_3 = \frac{\partial}{\partial \varphi}, 
\]
with $a = b = 0, c = - 3\varepsilon(\varphi + 1)\theta$. 
In fact we have exactly $[\zeta_1, \zeta_2] = \xi_2 + \varphi\xi_4 =: \zeta_4, 
[\zeta_1, \zeta_3] = 0, [\zeta_2, \zeta_3] = - \xi_3 - \theta\xi_4 =: \zeta_5, [\zeta_1, \zeta_4] = 0, 
[\zeta_1, \zeta_5] = - \xi_4 =: \zeta_6, [\zeta_2, \zeta_4] = 0, [\zeta_3, \zeta_4] = \xi_4, 
[\zeta_2, \zeta_5] = \xi_5 - 2\theta\xi_6 - \theta^2\xi_7 =: \zeta_7, 
[\zeta_1, \zeta_6] = 0, [\zeta_1, \zeta_7] = 2(- \xi_6 - \theta\xi_7), 
[\zeta_2, \zeta_6] = - \xi_6 - \theta_7 =: \zeta_8, 
[\zeta_3, \zeta_6] = 0, [\zeta_3, \zeta_7] = 0, [\zeta_2, \zeta_7] = -3\varepsilon(1+\varphi)(\zeta_8 + 
\theta\zeta_9) \equiv -3\varepsilon(1+\varphi)\theta \zeta_9$ ${\mbox {\rm mod}}. \ {\mathcal E}^{(4)}$. 

By Proposition \ref{necessary-condition}, we have ${\mbox {\it{SVC}}} \subseteq C$, 
where $C$ is the characteristic cone
$$
C = \{ u_1\xi_1 + u_2\xi_2 + u_3\xi_3 + u_4\xi_4 \mid u_1u_4 - u_2u_3 = 0\}. 
$$

If $\varepsilon = 0$, then the prolongation $(Z, E)$ of 
$(X, D)$ has growth $(3, 5, 7, 8, 9)$ everywhere and $(\R^7, D)$ is a hyperbolic $(4, 7)$-distribution 
of $C_3$-type everywhere and in that case we have ${\mbox {\it{SVC}}} = C$ by Theorem \ref{main-theorem}. 

If $\varepsilon \not= 0$, then $(X, D)$ is of mixed type, and 
the $(3, 5, 7, 8, 9)$-locus $\Sigma$ of the prolongation $(Z, E)$ of $(X, D)$ 
is a non-empty hypersurface in $Z$ given by $\theta = 0, \varphi = -1$. 
Then we see by Theorem \ref{main-theorem3}(2), that 
$$
{\mbox {\it{SVC}}} \subseteq \{ u_1 = u_3 = 0\} \cup \{ u_2 = u_4 = 0\} \cup \{ u_1 = u_3, u_2 = u_4\}. 
$$

The exact determination of the singular velocity cones for hyperbolic $(4, 7)$-distributions of mixed type 
is an open problem as far as the authors know. 


\begin{flushleft}
Goo ISHIKAWA
\\
e-mail: ishikawa@math.sci.hokudai.ac.jp
\end{flushleft}

\begin{flushleft}
Yoshinori MACHIDA,
\\
e-mail: machida.yoshinori@shizuoka.ac.jp
\end{flushleft}

\end{document}